\pdfoutput=1
\RequirePackage{ifpdf}
\ifpdf % We~are running pdfTeX in pdf mode
\documentclass[pdftex]{sigma}
\else
\documentclass{sigma}
\fi

\numberwithin{equation}{section}

\newtheorem{Theorem}{Theorem}[section]
\newtheorem{Corollary}[Theorem]{Corollary}
\newtheorem{Proposition}[Theorem]{Proposition}
\newtheorem{Question}[Theorem]{Question}
 { \theoremstyle{definition}
\newtheorem{Definition}[Theorem]{Definition}
\newtheorem{Remark}[Theorem]{Remark} }

\newcommand{\R}{\mathbb{R}}

\newcommand{\Ric}{\operatorname{Ric}}
\newcommand{\Rm}{\operatorname{Rm}}

\newcommand{\id}{\operatorname{id}}
\newcommand{\loc}{\text{loc}}

\newcommand{\diam}{\operatorname{diam}}
\newcommand{\met}{\operatorname{met}}
\newcommand{\G}{\operatorname{G}}
\newcommand{\Rf}{\operatorname{Rf}}
\renewcommand{\L}{\mathcal{L}}

\begin{document}
\allowdisplaybreaks

\newcommand{\arXivNumber}{2007.14967}

\renewcommand{\thefootnote}{}

\renewcommand{\PaperNumber}{128}

\FirstPageHeading

\ShortArticleName{Defining Pointwise Lower Scalar Curvature Bounds for $C^0$ Metrics}

\ArticleName{Defining Pointwise Lower Scalar Curvature Bounds\\ for $\boldsymbol{C^0}$ Metrics with Regularization by Ricci Flow\footnote{This paper is a~contribution to the Special Issue on Scalar and Ricci Curvature in honor of Misha Gromov on his 75th Birthday. The full collection is available at \href{https://www.emis.de/journals/SIGMA/Gromov.html}{https://www.emis.de/journals/SIGMA/Gromov.html}}}

\Author{Paula BURKHARDT-GUIM}

\AuthorNameForHeading{P.~Burkhardt-Guim}

\Address{Department of Mathematics, University of California, Berkeley, USA}
\Email{\href{mailto:paulab@math.berkeley.edu}{paulab@math.berkeley.edu}}
\URLaddress{\url{https://math.berkeley.edu/~paulab/}}

\ArticleDates{Received July 30, 2020, in final form November 19, 2020; Published online December 04, 2020}

\Abstract{We survey some recent work using Ricci flow to create a class of local definitions of weak lower scalar curvature bounds that is well defined for $C^0$ metrics. We discuss several properties of these definitions and explain some applications of this approach to questions regarding uniform convergence of metrics with scalar curvature bounded below. Finally, we consider the relationship between this approach and some other generalized notions of lower scalar curvature bounds.}


\Classification{53E20; 53C21}

\renewcommand{\thefootnote}{\arabic{footnote}}
\setcounter{footnote}{0}

\section{Introduction}\label{section1}

In \cite{Gro} Gromov presented the following theorem on the uniform limit of metrics satisfying a~pointwise lower bound on the scalar curvature; cf.~\cite[p.~1118]{Gro} and \cite[Theorem~1]{Bam}:
\begin{Theorem}\label{thm:C0limit}
Let $M$ be a smooth manifold and $\kappa\colon M\to R$ a continuous function on $M$.
Suppose $g_i$ is a sequence of $C^2$ metrics on $M$ that converges locally uniformly to a $C^2$ metric $g$ on $M$. If $R(g_i) \geq \kappa$ everywhere on $M$ for $i = 1,2,\dots$, then $R(g) \geq \kappa$ everywhere on $M$ as well.
\end{Theorem}
The proof described by Gromov involves formulating positive scalar curvature as a $C^0$ quantity, by considering the mean convexity and dihedral angles of small cubes. Bamler provided an alternative proof of Theorem~\ref{thm:C0limit} in~\cite{Bam}, which used the evolution of the scalar curvature under Ricci flow and some results of Koch and Lamm~\cite{KL1} concerning the Ricci--DeTurck flow for a class of possibly nonsmooth initial data on Euclidean space.

In light of Bamler's approach to Theorem \ref{thm:C0limit}, it is natural to ask whether it is possible to use Ricci flow to formulate a generalized definition of lower scalar curvature bounds for $C^0$ metrics. In~\cite{PBG} the author proposed a class of such definitions, and proved some related results. Generally speaking, the lower scalar curvature bounds in~\cite{PBG} are determined by using the Ricci flow to ``regularize'' the singular metric and then observing the scalar curvature of the flow at small positive times. The purpose of this survey is to describe more precisely how to formulate such a definition, to discuss several natural properties of this definition, and to explain how to prove some stability and rigidity results concerning the uniform convergence of metrics satisfying certain lower scalar curvature bounds, such as the following:
\begin{Theorem}\label{thm:C2kappaapproximation}
Let $g$ be a $C^0$ metric on a closed manifold $M$ which admits a uniform approximation by $C^2$ metrics $g^i$ that have $R(g^i)\geq \kappa_i$, where $\kappa_i$ is some sequence of numbers such that $\kappa_i\xrightarrow[i\to \infty]{} \kappa$ for some number $\kappa$. Then $g$ admits a uniform approximation by smooth metrics with scalar curvature bounded below by $\kappa$.
\end{Theorem}
\begin{Corollary}\label{cor:torusstability}
Let $g$ be a smooth metric on the torus which admits a uniform approximation by smooth metrics $g^i$ such that $R(g^i)\geq -\varepsilon_i\xrightarrow[i\to\infty]{} 0$. Then $g$ is flat.
\end{Corollary}

\section{Preliminaries: Ricci and Ricci--DeTurck flow}
If $M$ is a smooth manifold and $g_0$ is a smooth Riemannian metric on $M$, the Ricci flow starting from $g_0$ is a solution to
\begin{gather*}
\partial_t \bar g = -2\Ric(\bar g) \qquad \text{in}\quad M \times (0,T),\\
\bar g(0) = g_0,
\end{gather*}
where $\bar g(t)$ is a smooth, time-dependent family of Riemannian metrics on a $M$ for $t\in (0,T)$. If~$M$ is closed it is known that a short-time solution to the Ricci flow always exists and is unique; see \cite[Theorems~5.2.1 and~5.2.2]{Top}.

We will also find it useful to consider the Ricci--DeTurck flow, introduced by DeTurck in \cite{De}, a strongly parabolic flow that is related to the Ricci flow by pullback via a family of diffeomorphisms. More specifically, we define the following operator, which maps symmetric $2$-forms on~$M$ to vector fields:
\begin{equation*}%\label{eq:Xoperator}
X_{\bar g}(g):= \sum_{i=1}^n\big(\nabla^{\bar g}_{e_i}e_i - \nabla^{g}_{e_i}e_i\big),
\end{equation*}
where $\{e_i\}_{i=1}^n$ is any local orthonormal frame with respect to $g$. Then the Ricci--DeTurck equation is
\begin{equation}\label{eq:RDTeq}
\partial_t g(t) = -2\Ric(g(t)) - \L_{X_{\bar g(t)}(g(t))}g(t),
\end{equation}
where $\bar g(t)$ is a background Ricci flow. As mentioned, if $g(t)$ solves (\ref{eq:RDTeq}) then it is related to a~Ricci flow via pullback by diffeomorphisms. More precisely, if $g(t)$ solves (\ref{eq:RDTeq}) and $\chi_t\colon M\to M$ is a family of diffemorphisms satisfying
\begin{gather}
X_{\bar g(t)}(g(t))f = \frac{\partial}{\partial t}(f\circ\chi_t) \text{ for all } f\in C^\infty(M),\nonumber\\
\chi_0 = \id,\label{eq:diffeoseq}
\end{gather}
 then $\chi_t^*g(t)$ solves the Ricci flow equation with initial condition $g(0)$.

Let $g_t$ be a solution to the Ricci--DeTurck equation, and write $g_t = \bar g_t + h_t$, where $\bar g_t$ is the smooth background Ricci flow from (\ref{eq:RDTeq}). Then the evolution equation (see \cite[Appendix~A]{BK} for more details, and a more general setting) for $h_t$ is
\begin{equation}\label{eq:hevolution}
\partial_t h_t + Lh_t = Q[h_t],
\end{equation}
where the linear part, $L$, is
\begin{equation*}\label{eq:Lis}
L h_t := \Delta^{\bar g_t}h_t + 2\Rm^{\bar g_t}(h_t) := \Delta^{\bar g_t}h_t + 2{\bar g}^{pq}R_{pij}^{m}h_{q m}{\rm d}x^i\otimes {\rm d}x^j,
\end{equation*}
and the quadratic term $Q$ may be written as
\begin{equation*}
Q[h_t] = Q^0_t + \nabla^* Q^1_t,
\end{equation*}
where
\begin{equation*}%\label{eq:Q0is}
Q^0_t = (\bar g + h)^{-1}\star (\bar g + h)^{-1} \star \nabla h \star \nabla h + \big((\bar g + h)^{-1} - \bar g^{-1}\big)\star \Rm^{\bar g_t}\star h
\end{equation*}
and
\begin{equation*}%\label{eq:Q1is}
\nabla^*Q^1_t:= \nabla_p \big(\big((\bar g+h)^{pq} - \bar g^{pq}\big)\nabla_qh_{ij}\big) = \nabla\big(\big((\bar g + h)^{-1} - \bar g^{-1}\big)\star \nabla h\big).
\end{equation*}
Here we use the notation $A\star B$ for two tensor fields $A$ and $B$ to mean a linear combination of products of the coefficients of $A$ and $B$, and $(\bar g + h)^{-1}$ and $\bar g^{-1}$ denote tensor fields with coefficients $(\bar g + h)^{ij}$ and $\bar g^{ij}$ respectively. We sometimes also write $h_t$ using the integral representation
\begin{gather}
h_t(x) = \int_M \bar K(x,t;y, 0)h_0(y){\rm d}\bar g_0(y)\nonumber\\
\hphantom{h_t(x) =}{} + \int_{M\times[0,t]}\big(\bar K(x,t; y, s)Q^0_s(y)
+ \nabla^* \bar K(x,t;y,s)Q^1_s(y)\big){\rm d}\bar g_s(y),\label{eq:integraleq}
\end{gather}
where $\bar K$ denotes the heat kernel for the background Ricci flow $\bar g(t)$. We refer to~(\ref{eq:hevolution}) as the Ricci--DeTurck perturbation equation.

There are several results on the existence of Ricci--DeTurck flows starting from $C^0$ initial data, such as \cite[Theorem $1.1$]{Sim}, \cite[Theorem $4.3$]{KL1}, and \cite[Theorem $5.3$]{KL2}. Because we will consider flows from initial metrics that are close to one another in the $C^0$ sense, the result that is most relevant to this survey is (cf.\ \cite[Lemma $3.2$ and Corollary $3.3$]{PBG}):

 \begin{Proposition}\label{prop:RDTexists}
Let $M^n$ be a smooth, closed manifold. There exist constants $\varepsilon = \varepsilon(n)$ and $C = C(n)$ such that the following is true:

For every metric $g_0\in C^0(M)$ and every smooth background metric $\bar g_0$, if $||g_0 - \bar g_0||_{L^{\infty}(M)} < \varepsilon$ and $\bar g_t$ is the Ricci flow starting from $\bar g_0$, then there exists $T = T(\bar g(t))$ sufficiently small so that there is a solution $g_t$ to the integral equation~\eqref{eq:integraleq} such that
\begin{equation*}
||g_t - \bar g_t||_{L^\infty(M\times (0,T))} \leq C||g_0 - \bar g_0||_{L^{\infty}(M)}.
\end{equation*}

Moreover, by taking $\varepsilon$ and $T$ smaller, we may find constants $c_k$ depending only on $k$, the dimension, and bounds for the derivatives of $\Rm(\bar g)$, such that the solution $g_t$ is smooth on $M\times (0,T]$, continuous on $M\times [0,T]$, and satisfies
\begin{equation*}
\big|\nabla^{k}(g_t - \bar g_t)\big| \leq \frac{c_k}{t^{k/2}}||g_0 - \bar g_0||_{L^{\infty}(M)}
\end{equation*}
for all $t\in (0,T']$, where $\nabla$ denotes the covariant derivative with respect to $\bar g(t)$.
 \end{Proposition}
 The existence of a solution to (\ref{eq:integraleq}) is proven in a similar fashion to \cite{KL1}: in~\cite{KL1} Koch and Lamm construct suitable Banach spaces so that a solution to~(\ref{eq:integraleq}) arises from an application of the Banach fixed point theorem. In~\cite{PBG} one performs a similar construction, working on a~Ricci flow background rather than a stationary background. Once the existence of a solution to the integral equation~(\ref{eq:integraleq}) has been established, smoothness of the solution and bounds on the derivatives may be proven by iterative application of parabolic interior estimates.

Under the Ricci flow, the scalar curvature evolves by
\begin{equation}\label{eq:Revolution}
\partial_t R = \Delta^{\bar g(t)} R + 2|\Ric|^2;
\end{equation}
see \cite[Proposition~2.5.4]{Top}.
Making an orthogonal decomposition, we may conclude that
\begin{equation}\label{eq:diffineqR}
\partial_t R \geq \Delta^{\bar g(t)}R + \frac{2}{n}R^2;
\end{equation}
this is \cite[Corollary~2.5.5]{Top}. For smooth Ricci flows on a closed manifold, a general bound that follows from the maximum principle is (see \cite[Corollary~3.2.5]{Top})
\begin{equation*}%\label{eq:classicalRlowerbound}
R(x,t) \geq -\frac{n}{2t},
\end{equation*}
for all $(x,t)\in M\times (0,T]$.

If, instead, $g(t)$ is a Ricci--DeTurck flow, then recall that $g(t) = (\chi_t^{-1})^*\bar g(t)$ for some Ricci flow $\bar g(t)$, where the family $(\chi_t)$ satisfies (\ref{eq:diffeoseq}) and $X$ is the corresponding vector field. Therefore, pushing forward (\ref{eq:diffineqR}) by $(\chi_t)$ we find that, under the Ricci--DeTurck flow,
\begin{equation*}%\label{eq:RDTRevolution}
\partial_t R \geq \Delta^{g(t)}R - \langle X, \nabla R\rangle + \frac{2}{n}R^2;
\end{equation*}
see also \cite[p.~6]{Bam}. It follows from the maximum principle that, if $g(t)$ is a Ricci or Ricci--DeTurck flow on a closed manifold, starting from a smooth initial metric, and defined on the interval~$[0,T]$, then if $R^{g_0} \geq \kappa_0\in \R$, we have (cf.\ \cite[Theorem~3.2.1]{Top})
\begin{equation}\label{eq:Rpreservation}
R^{g(t)} \geq \frac{\kappa_0}{1 - \big(\frac{2\kappa_0}{n}t\big)} \geq \kappa_0,
 \end{equation}
 for all $t\in [0,T]$.

\section[Regularizing Ricci flow and pointwise lower scalar curvature bounds]{Regularizing Ricci flow\\ and pointwise lower scalar curvature bounds}

As mentioned in Section~\ref{section1}, in order to determine the generalized lower scalar curvature bound for a~$C^0$ metric we first regularize the metric by the Ricci flow and then study the scalar curvature of the flow for small positive times. It is not obvious what is meant by a Ricci flow starting from~$C^0$ initial data; however, we have \cite[Theorem~1.1]{PBG}:

\begin{Theorem}\label{thm:RRF}
Let $M$ be a closed manifold and $g_0$ a $C^0$ metric on $M$. Then there exists a~time-dependent family of smooth metrics $(\tilde g_t)_{t\in(0,T]}$ and a continuous surjection $\chi\colon M\to M$ such that the following are true:
\begin{enumerate}\itemsep=0pt
\item[$(a)$] The family $(\tilde g_t)_{t\in(0,T]}$ is a Ricci flow, and
\item[$(b)$] There exists a smooth family of diffeomorphisms $(\chi_t)_{t\in (0,T]}: M\to M$ such that
\begin{equation*}
\chi_t\xrightarrow[t\to 0]{C^0}\chi\qquad \text{and} \qquad ||(\chi_t)_*\tilde g_t - g_0||_{C^0(M)}\xrightarrow[t\to 0]{}0.
\end{equation*}
\end{enumerate}
Moreover, for any $x\in M$, $\diam \{\chi_s(x)\colon s\in(0,t]\} \leq C\sqrt{t}$ for some constant $C>0$ independent of $x$, where the diameter is measured with respect to a fixed smooth background metric, and any two such families are isometric, in the sense that if $\tilde g_t'$ is another such family with corresponding continuous surjection $\chi'$, then there exists a stationary diffeomorphism $\alpha\colon M\to M$ such that $\alpha^* \tilde g_t = \tilde g_t'$ and $\chi\circ\alpha = \chi'$.
\end{Theorem}
The pair $((\tilde g_t)_{t\in(0,T]}, \chi)$ is called a \emph{regularizing Ricci flow} for $g_0$.
\begin{Remark}
It has not been shown that there are metrics for which $\chi$ cannot be taken to be a homeomorphism, though we suspect that such metrics do exist.
\end{Remark}

\begin{Remark}\label{rmk:broaderuniqueness}
Moreover, the regularizing Ricci flow is invariant under $C^0$ isometry in a more general sense: if $\varphi\colon \big(M_2, g^2\big) \to \big(M_1, g^1\big)$ is a metric space isometry, then, for any two regularizing Ricci flows $\big(\big(\tilde g^1(t)\big)_{t\in (0,T^1]}, \chi^1\big)$ and $\big(\big(\tilde g^2(t)\big)_{t\in (0,T^2]}, \chi^2\big)$ for $g^1$ and $g^2$ respectively, there is a~stationary diffeomorphism $\alpha\colon M_2 \to M_1$ such that $\alpha^*\tilde g^1(t) = \tilde g^2(t)$ for all $t\in \big(0,\min\big\{T^1,T^2\big\}\big]$ and $\chi^1\circ\alpha = \varphi \circ \chi^2$; this is \cite[Corollary~5.5]{PBG}.
\end{Remark}

In the smooth case, we may take $(\tilde g_t)$ to be the classical Ricci flow by choosing $(\chi_t)$ to be the family of diffeomorphisms which solves~(\ref{eq:diffeoseq}) and is defined for $t\in [0,T]$, and then pulling back the Ricci--DeTurck flow by $(\chi_t)_{t\in [0,T]}$. In the case of $C^0$ initial data, it is still possible to solve the differential equation in (\ref{eq:diffeoseq}) subject to the ``initial condition'' $\chi_{t_0} = \id$ for some $t_0> 0$, since~$X$ is nonsingular for $t_0 > 0$, but it is not possible to prescribe $\chi_0 = \id$, since $X$ may be singular at time~$0$. Instead, we take $\chi$ to be the uniform limit as $t\searrow 0$ of the solution $(\chi_t)_{t\in (0,T]}$ to the differential equation in~(\ref{eq:diffeoseq}) with $\chi_{t_0} = \id$. The family $(\tilde g_t)$ is then constructed by pulling back the Ricci--DeTurck flow in Proposition~\ref{prop:RDTexists} by $(\chi_t)_{t\in (0,T]}$.

The non-uniqueness of the regularizing Ricci flow is an artifact of the non-unique choice of ``initial data'' $\chi_{t_0}$. That the regularizing Ricci flow is unique up to isometry is due to the uniqueness of Ricci flows from smooth initial data on closed manifolds, and compactness of the isometry group for smooth Riemannian metrics on closed manifolds.

We now explain how to use the regularizing Ricci flow to define pointwise lower scalar curvature bounds for $C^0$ metrics. Henceforth we will refer to this definition as the ``weak sense''. A satisfactory definition of the weak sense should satisfy the following requirements: For any constant $\kappa$, we should have
\begin{enumerate}\itemsep=0pt
\item \textit{Stability under greater-than-second-order perturbation:} If $g'$ and $g''$ are two $C^0$ metrics that agree to greater than second order around a point $x_0$, i.e., if, for some fixed smooth background metric, we have $|g'(x) - g''(x)| \leq cd^{2+\eta}(x,x_0)$ for some $c, \eta>0$ and all $x$ in a neighborhood of $x_0$, then $g'$ should have scalar curvature bounded below by $\kappa$ in the weak sense at $x_0$ if and only if $g''$ does. Moreover, if $g'$ and $g''$ are $C^0$ metrics on different manifolds which merely agree to greater than second order under pullback by a locally defined diffeomorphism, the conclusion should still hold.
\item \textit{Preservation of global lower bounds under the Ricci flow:} If $g$ is a $C^0$ metric on a closed manifold that has scalar curvature bounded below by $\kappa$ in the weak sense at every point, and $\tilde g_t$ is a regularizing Ricci flow for $g$, then $\tilde g_t$ should have scalar curvature bounded below by $\kappa$ at every point for all $t>0$ for which the flow is defined. This is true for Ricci flows starting from smooth initial data.
\item \textit{Agreement with the classical notion for $C^2$ metrics:} If $g$ is a $C^2$ metric with scalar curvature bounded below by $\kappa$ at $x_0$ in the weak sense for $C^0$ metrics, then $g$ should have scalar curvature bounded below by~$\kappa$ at~$x_0$ in the classical sense. Conversely, if $g$ has scalar curvature bounded below by~$\kappa$ at~$x_0$ in the classical sense, then the same should hold in the weak sense.
\end{enumerate}

Perhaps the most natural approach would be to declare that $R(g_0)\geq \kappa$ at $x$ if $\lim\limits_{t\searrow 0}R(g_t, y) \geq \kappa$ for some $y\in \chi^{-1}(x)$ (we consider the scalar curvature at $y \in \chi^{-1}(x)$ rather than at $x$ itself because $g_t$ only converges to $g_0$ modulo pushforward by the family $(\chi_t)$). As we shall discuss below, we wish to leverage exponential heat kernel estimates in conjunction with (\ref{eq:Revolution}) to show that ($2$) holds, so we modify this approach by instead requiring that the lower bound holds in the limit, in a small time-dependent neighborhood of $y$, where the size of the neighborhood tends to $0$ as $t\searrow 0$:

\begin{Definition}\label{def:RRFscalar}
Let $M^n$ be a closed manifold and $g_0$ a $C^0$ metric on $M$. For $0 < \beta < 1/2$ we say that $g_0$ has scalar curvature bounded below by $\kappa$ at $x$ in the \emph{$\beta$-weak sense} if there exists a~regularizing Ricci flow $((\tilde g_t)_{t\in(0,T]}, \chi)$ for $g_0$ such that, for some point $y\in M$ with $\chi(y) = x$, we have
\begin{equation}\label{eqn:RRFscalar}
\inf_{C>0}\Big(\liminf_{t\searrow 0}\Big(\inf_{B_{\tilde g(t)}(y,Ct^\beta)}R^{\tilde g}(\cdot, t)\Big)\Big) \geq \kappa,
\end{equation}
where $B_{\tilde g(t)}\big(y,Ct^\beta\big)$ denotes the ball of radius $Ct^\beta$ about $y$, measured with respect to the metric~$\tilde g(t)$, and $R^{\tilde g}(\cdot,t)$ denotes the scalar curvature of $\tilde g_t$ .
\end{Definition}
\begin{Remark}\label{rmk:welldefinednessofR}
In fact, Definition~\ref{def:RRFscalar} is independent of choice of $y$, so it is equivalent to require that (\ref{eqn:RRFscalar}) hold at $y$ for \emph{all} $y$ with $\chi(y) = x$. Moreover, the fact that any two regularizing Ricci flows for $g_0$ are isometric (see Theorem \ref{thm:RRF}) implies that Definition~\ref{def:RRFscalar} holds for some regularizing Ricci flow if and only if it holds for all regularizing Ricci flows for~$g_0$.
\end{Remark}
\begin{Remark}
Definition \ref{def:RRFscalar} may also be formulated in terms of Ricci--DeTurck flow (see \cite[Lemma~6.4]{PBG}). The equivalent statement using Ricci--DeTurck flow is that there is a Ricci--DeTurck flow $g_t$ starting from $g_0$ in the sense of Proposition \ref{prop:RDTexists} and $\bar g_0$ a stationary metric on~$M$ that is uniformly bilipschitz to $(g_t)_{t\in (0,T]}$ such that
\begin{equation*}
\inf_{C>0}\Big(\liminf_{t\searrow 0}\Big(\inf_{B_{\bar g_0}(x, Ct^\beta)}R^{g_t}(\cdot)\Big)\Big) \geq \kappa.
\end{equation*}
\end{Remark}

We take the infimum over all $C>0$ to ensure that the condition that $R(g_0)\geq 0$ at $x$ in the $\beta$-weak sense is scaling invariant, since the condition $0<\beta<1/2$ disrupts the parabolic scaling of the expression. The significance of $\beta$ in this range is that away from $B_{\tilde g(t)}\big(y,Ct^\beta\big)$ the heat kernel for the flow decays exponentially as $t\searrow 0$; we will discuss how this is used below. It is unknown whether $\beta$ could instead be replaced by~$1/2$, or even to replace~(\ref{eqn:RRFscalar}) by $\liminf\limits_{t\searrow 0}R(g_t, y)$. It is clear that Definition~\ref{def:RRFscalar} satisfies item~($3$), since if $g_0$ is $C^2$, then by Remark~\ref{rmk:welldefinednessofR} we may take the regularizing Ricci flow as the usual Ricci flow with $\chi = \id$, and
\begin{equation*}
\inf_{C>0}\Big(\liminf_{t\searrow 0}\Big(\inf_{B_{\tilde g(t)}(y,Ct^\beta)}R^{\tilde g}(\cdot, t)\Big)\Big) = \lim_{t\to 0}R^{\tilde g}(x, t) = R^{g}(x).
\end{equation*}

Item ($1$) holds due to the following stability result for the scalar curvatures of regularizing Ricci flows from metrics which agree to greater than second order at a point; this is \cite[Theorem~1.4]{PBG}:
\begin{Theorem}\label{thm:secondorderagreement}
Suppose $g'$ and $g''$ are two $C^0$ metrics on closed manifolds $M'$ and $M''$ respectively, and that there is a locally defined diffeomorphism $\phi\colon U\to V$ where $U$ is a neighborhood of $x_0'$ in $M'$ and $V$ is a neighborhood of~$x_0''$ in~$M''$ with $\phi(x_0') = x_0''$. Suppose furthermore that $g'$ and $\phi^*g''$ agree to greater than second order around $x_0'$, i.e., with respect to some fixed smooth background metric, $|g'(x) - \phi^*g''(x)| \leq cd^{2+\eta}(x,x_0)$ for some $c, \eta>0$ and all $x$ in a neighborhood of~$x_0'$. Then there exist regularizing Ricci flows $(\tilde g_t', \chi')$ and $(\tilde g_t'', \chi'')$ for $g'$ and $g''$ respectively such that, for $1/(2+\eta) < \beta <1/2$, $C>0$, and~$t$ sufficiently small depending on $C$, $\beta$, and $\eta$, we have
\begin{equation*}
\sup_{B(x_0', Ct^\beta)}\big|R^{(\chi_t')_*\tilde g'_t} - \phi^*R^{(\chi_t'')_*\tilde g''_t}\big|\leq ct^\omega,
\end{equation*}
where $\omega$ is some positive exponent, $c$ is a constant that does not depend on $t$ or $C$, $R^{(\chi_t')_*\tilde g'_t}$ and $R^{(\chi_t'')_*\tilde g''_t}$ denote the scalar curvatures with respect to $(\chi_t')_*\tilde g_t'$ and $(\chi_t'')_*\tilde g_t''$ respectively, and $(\chi_t')$ and $(\chi_t'')$ are the smooth families of diffeomorphisms for $\tilde g_t'$ and $\tilde g_t''$ respectively, whose existence is given by~(b) in Theorem~{\rm \ref{thm:RRF}}.

In particular, Definition~{\rm \ref{def:RRFscalar}} holds for $g'$ at $x_0'$ if and only if it holds for $g''$ and $x_0''$.
\end{Theorem}
\begin{Remark}
By Theorem \ref{thm:secondorderagreement}, Definition~\ref{def:RRFscalar} descends to the space of germs of metrics at a~point, and further descends to the quotient space induced on the space of germs of metrics at $x_0$ by the relation $[g] \sim [g']$ if $g$ and $g'$ agree to greater than second order at~$x_0$. One may also use this fact to define weak pointwise lower scalar curvature bounds for $C^0$ metrics on open manifolds.
\end{Remark}

The strategy for proving Theorem \ref{thm:secondorderagreement} is essentially to endow the Banach spaces defined by Koch and Lamm in \cite{KL1} with the weight $\omega(x,t) = \max\big\{\big(d(x_0, x) + \sqrt{t}\big)^{-2-\eta}, 1\big\}$ to offset the evolution of two metrics that initially agree to greater than second order about $x_0$, and then to show that $||\omega(\cdot, t)(g_t' - \phi^*g_t'')||_{C^0(M)} \leq C(||\omega(\cdot, 0)g_0' - \phi^*g_0''||) < \infty$. Closeness of the scalar curvatures is obtained from the $C^0$ estimate by interpolating between higher and lower order derivatives.

Moreover, we have that item ($2$) holds; this is \cite[Theorem~1.5]{PBG}:
\begin{Theorem}\label{thm:preservation}
Suppose that $g_0$ is a $C^0$ metric on a closed manifold $M$, and suppose there is some $\beta \in (0,1/2)$ such that $g_0$ has scalar curvature bounded below by $\kappa$ in the $\beta$-weak sense at all points in~$M$. Suppose also that $(\tilde g(t))_{t\in (0,T]}$ is a Ricci flow starting from $g_0$ in the sense of Theorem~{\rm \ref{thm:RRF}}. Then the scalar curvature of $\tilde g(t)$, $R(\tilde g(t))$, satisfies $R(\tilde g(t)) \geq \kappa$ everywhere on~$M$, for all $t\in (0,T]$.
\end{Theorem}

To prove Theorem \ref{thm:preservation} one supposes that the theorem is false and then iteratively uses the heat kernel estimates and (\ref{eq:Revolution}) to show that there must be a sequence $(x_i, t_i)$ with $t_i\searrow 0$ and all $x_i\in B_{\tilde g(t_i)}\big(x_0, t_i^{\beta}\big)$ for some $x_0$, such that $R(x_i, t_i)$ is bounded above away from $\kappa$ for all~$i$, so that Definition~\ref{def:RRFscalar} must fail at~$x_0$.

\section{Some applications: global lower bounds}

Recall that the purpose of Definition \ref{def:RRFscalar} was to formulate a \emph{pointwise} lower scalar curvature bound for a $C^0$ metric $g$. Nevertheless, if $g$ satisfies a (constant) global lower bound on the scalar curvature in the sense of Definition~\ref{def:RRFscalar}, Theorem~\ref{thm:preservation} leads to several alternate characterizations of $g$, which we will discuss in this section.

 We have (cf.\ \cite[Corollary~1.6]{PBG}):
 \begin{Theorem}\label{thm:globalequivalence}
 Let $g$ be a $C^0$ metric on a closed manifold $M$. Let $\kappa$ be some constant. The following are equivalent:
 \begin{enumerate}\itemsep=0pt
 \item[$1.$] There exists some $\beta < 1/2$ such that $R(g)\geq \kappa$ in the $\beta$-weak sense everywhere on $M$.
 \item[$2.$] For any regularizing Ricci flow $(\tilde g(t))_{t\in (0,T]}$ for $g$, $R(\tilde g(t))\geq \kappa$ everywhere on $M$, for all $t\in (0,T]$.
 \item[$3.$] There exists a sequence of smooth metrics $g^i$ on $M$ such that $R(g^i)\geq \kappa$ everywhere on $M$ for $i = 1,2, \dots$, and $g^i\xrightarrow[i\to\infty]{C^0} g$.
\end{enumerate}
 \end{Theorem}
 \begin{Remark}
In particular, Theorem \ref{thm:globalequivalence} implies that if the global lower bound $R(g)\geq \kappa$ holds in the $\beta$-weak sense for some value of $\beta < 1/2$, then it holds for all $\beta < 1/2$.
\end{Remark}

 To see why Theorem \ref{thm:globalequivalence} is true, first choose a sequence of times $t_i\searrow 0$, and let $\tilde g(t)$ be a regularizing Ricci flow for $g$. Then Theorem \ref{thm:RRF} says that $(\chi_{t_i})_*\tilde g(t_i)$ converges uniformly to $g$. Moreover, Theorem \ref{thm:preservation} guarantees that $R(\tilde g(t_i)) \geq \kappa$ for all $i$, so the same is true for $R((\chi_{t_i})_*\tilde g(t_i))$, so $(1)\Rightarrow (3)$.

Furthermore, if there exists a sequence of smooth approximators $g^i$ as in $(3)$, then for sufficiently large $i$ all $g^i$ are close enough to some smooth background metric $\bar g_0$ so that there exists a uniform value $T$ such that the Ricci--DeTurck flows $g^i(t)$ starting from $g^i$ in the sense of Proposition~\ref{prop:RDTexists} are defined for $t\in (0,T]$. Moreover, $R(g^i(t))\geq \kappa$ for all $t\in (0,T]$ and $i$, according to~(\ref{eq:Rpreservation}).

By the derivative estimates given by Proposition~\ref{prop:RDTexists}, the flows $g^i(t)$ converge smoothly to the Ricci--DeTurck flow $g(t)$ for $g$ on $M$ for all fixed $t>0$. In particular, for any fixed $t\in (0,T]$ and all $i$ sufficiently large, we have $R(g(t)) \geq \kappa$. Since $g(t)$ is related to a regularizing Ricci flow $\tilde g(t)$ for $g$ by a family of diffeomorphisms, we find that $R(\tilde g(t))\geq \kappa$ for all $t\in (0,T]$, so $(3)\Rightarrow (2)$. It is clear that $(2)\Rightarrow (1)$ for any $\beta < 1/2$. In particular, we have shown that $(1)\Rightarrow (3) \Rightarrow (2) \Rightarrow (1)$, which proves Theorem \ref{thm:globalequivalence}.

 Observe that in our argument to show that $(3)\!\Rightarrow\! (2)$ it was only necessary to have \mbox{$R(g^i(t)) \!\geq\! \kappa_i$} where $\kappa_i \to \kappa$, since by smooth convergence of the time slices we would then have $R(g(t))\geq \kappa_i$ for all $i$ sufficiently large, and hence $R(g(t)) \geq \kappa$. In particular, Theorems~\ref{thm:globalequivalence} and~\ref{thm:preservation} imply that it is sufficient to require that the approximating metrics $g^i$ be only $C^0$ with $R(g^i)\geq \kappa_i$ in the sense of Definition~\ref{def:RRFscalar}. In particular, we have the following ``continuity property'' for global lower bounds in the sense of Definition~\ref{def:RRFscalar} (cf.~\cite[Theorem~1.7]{PBG}):
 \begin{Theorem}\label{thm:kappaapproximation}
Let $g$ be a $C^0$ metric on a closed manifold $M$ which admits a uniform approximation by $C^0$ metrics $g^i$ such that there is some $\beta< 1/2$ so that $g^i$ has scalar curvature bounded below by $\kappa_i$ in the $\beta$-weak sense everywhere on $M$, where $\kappa_i$ is some sequence of numbers such that $\kappa_i\xrightarrow[i\to \infty]{} \kappa$ for some number $\kappa$. Then $g$ has scalar curvature bounded below by $\kappa$ in the $\beta$-weak sense. In particular, any regularizing Ricci flow $(\tilde g(t))_{t\in (0,T]}$ for $g$ satisfies $R(\tilde g(t)) \geq \kappa$ for all $t\in (0,T]$, so $g$ admits a uniform approximation by smooth metrics with scalar curvature bounded below by $\kappa$.
\end{Theorem}

Theorem \ref{thm:C2kappaapproximation} follows from Theorem \ref{thm:kappaapproximation} in the case where the $g^i$ are $C^2$. We remark that Theorem~\ref{thm:C2kappaapproximation} may also be proved directly by using the Ricci--DeTurck flow as in $(3)\Rightarrow (2)$ in the proof of Theorem~\ref{thm:globalequivalence}, using Proposition~\ref{prop:RDTexists} or~\cite[Theorem~1.1]{Sim}. By setting $\kappa_i = -\varepsilon_i$, Theorem~\ref{thm:C2kappaapproximation} answers in the affirmative the following question, posed by Gromov in \cite{Gro}:
\begin{Question}[{\cite[p.~1119]{Gro}}]\label{q:Gromov}
Let $g$ be a continuous Riemannian metric on a closed manifold~$M$ which admits a $C^0$-approximation by smooth Riemannian metrics $g_i$ with $R(g_i)\geq -\varepsilon_i\xrightarrow[i\to \infty]{}0$. Does~$M$ admit a smooth metric of nonnegative scalar curvature?
\end{Question}
Moreover, Theorem \ref{thm:C2kappaapproximation} implies that $g$ admits a uniform approximation by smooth metrics with nonnegative scalar curvature.

In light of the scalar torus rigidity theorem, which says that any Riemannian manifold with nonnegative scalar curvature that is diffeomorphic to the torus must be isometric to the flat torus (see~\cite{GL,SY}), it is natural to ask whether Definition~\ref{def:RRFscalar} satisfies an analogous rigidity result for~$C^0$ metrics on the torus. If $g$ is a $C^0$ metric on the torus with nonnegative scalar curvature in the $\beta$-weak sense, for some value of $\beta < 1/2$, and $(\tilde g_t, \chi_t)$ is a regularizing Ricci flow for~$g$, then Theorem~\ref{thm:preservation} implies that $R(\tilde g(t))\geq 0$ for all $t>0$, and hence the usual scalar torus rigidity theorem implies that $\tilde g(t)$ is flat for all $t>0$. In particular, the Ricci flow equation implies that $\tilde g(t)\equiv \tilde g$ for some stationary smooth flat metric $\tilde g$. Then $(\chi_t)_{\tilde g}\xrightarrow[t\to 0]{C^0} g$ implies that~$\tilde g$ is equivalent to~$g$ in the Gromov--Hausdorff sense. We have shown:
\begin{Corollary}\label{cor:torusrigidity}
Suppose $g$ is a $C^0$ metric on the torus $\mathbb{T}$, and that there is some $\beta \in (0,1/2)$ such that $g$ has nonnegative scalar curvature in the $\beta$-weak sense everywhere. Then $(\mathbb{T}, g)$ is isometric as a metric space to the standard flat metric on $\mathbb{T}$.
\end{Corollary}
Corollary~\ref{cor:torusrigidity} and Theorem~\ref{thm:C2kappaapproximation} together imply Corollary~\ref{cor:torusstability}.

\begin{Remark}
Corollary \ref{cor:torusrigidity} is in fact the optimal result, i.e., it is not possible to show that there is a \emph{Riemannian} isometry between $g$ and the standard flat metric. In the case where $g^1$ and $g^2$ are smooth metrics, a metric space isometry is automatically a smooth Riemannian isometry. However, there exist examples of $C^0$ isometries between $C^0$ metrics which are not $C^1$. Moreover, by Remark \ref{rmk:broaderuniqueness}, the regularizing Ricci flow is invariant under $C^0$ isometry, so it is not possible to use regularizing Ricci flow to distinguish Riemannian isometries from isometries which are merely metric space isometries.
\end{Remark}

\section{Relation to some notions of positive scalar curvature}

We will now discuss the relationship between Definition \ref{def:RRFscalar} and some other possible definitions of lower scalar curvature bounds for $C^0$ metrics, both in the global and local setting. Throughout we take $M$ to be a closed manifold. We first define the following classes of $C^0$ metrics on $M$, each of which is a class of $C^0$ metrics that have scalar curvature globally bounded below by a~constant~$\kappa$ in some reasonable generalized sense:
\begin{gather*}
C^0_{\beta}(M,\kappa) := \big\{g\in C^0(M)\colon R(g)\geq \kappa \text{ in the $\beta$-weak sense}\big\},\\
C^0_{\met}(M,\kappa) := \Big\{g\in C^0(M)\colon \text{ there exists } g^i\in C^2 \text{ with } R(g^i)\geq \kappa \text{ and } g^i\xrightarrow[i\to\infty]{C^0} g\Big\},\\
C^0_{\Rf}(M,\kappa) := \big\{g\in C^0(M) \colon \exists\, (\tilde g(t))_{t\in (0,T]} \text{ a regularizing Ricci flow for } g\\
\hphantom{C^0_{\Rf}(M,\kappa) := \big\{}{}\text{with } R(\tilde g(t)) \geq \kappa\, \forall\, t\in (0,T]\big\},\\
C^0_{\G}(M, \kappa) := \big\{c\in C^0(M) \colon R(g)\geq \kappa \text{ in the sense of Gromov \cite{Gro}}\big\}.
\end{gather*}

It is natural to study $C^0_{\Rf}(M, \kappa)$ in light of \cite{Bam}.
A metric $g$ is in $C^0_{\G}(M,0)$ if essentially it does not contain a cube with strictly mean convex faces, such that the dihedral angles are acute, or more generally, it is in $C^0_{\G}(M, \kappa)$ if the product metric of its product with an appropriate space form $S_{\kappa}$ is in $C^0_{\G}(M\times S_{-\kappa}, 0)$; see \cite[p.~1119]{Gro}. This is a natural definition because, if such a~cube were to exist in a manifold with (classical) nonnegative scalar curvature, then, as in \cite[pp.~1144--1145]{Gro}, one could glue together $2n$ copies of the cube and obtain a non-flat torus of nonnegative scalar curvature, contradicting the scalar torus rigidity theorem (\cite[Corollary~2]{SY} and \cite[Corollary~A]{GL}).

By Theorem \ref{thm:globalequivalence} we have that for any $\beta < 1/2$, $C^0_{\beta}(M,\kappa) = C^0_{\met}(M, \kappa) = C^0_{\Rf}(M,\kappa)$. The Ricci flow proof of Theorem~\ref{thm:C0limit} suggests a relationship between $C^0_{\beta}(M, \kappa) = C^0_{\Rf}(M,\kappa)$ and~$C^0_{\G}(M,\kappa)$. We have that $C^0_{\met}(M,\kappa) \subset C^0_{\G}(M,\kappa)$, since $C^0_{\G}(M,\kappa)$ contains all~$C^2$ metrics~$g$ with $R(g)\geq \kappa$ and is closed in $C^0$. Thus, $C^0_{\Rf}(M,\kappa) \subset C^0_{\G}(M,\kappa)$. One key question is whether, if a $C^0$ metric on a closed manifold has scalar curvature bounded below in the sense of \cite{Gro}, it necessarily has scalar curvature bounded below under the Ricci flow:

\begin{Question}
Suppose that $M$ is closed. Is $C^0_{\G}(M,\kappa)\subset C^0_{\Rf}(M,\kappa)$?
\end{Question}

If so, then this would imply that all $C^0$ notions of a global lower scalar curvature bound that we have mentioned would agree.

We remark that, aside from the classes of $C^0$ metrics listed above, Lee and LeFloch \cite{LL} have introduced a notion of positive distributional scalar curvature for metrics in $C^0\cap W^{1,n}_{\loc}$, which they have used to prove a version of the Positive Mass Theorem. In this paper we work with~$C^0$ metrics without any assumptions regarding their distributional derivatives, but in the case of metrics in $C^0\cap W^{1,n}_{\loc}$, it is an open question whether the definition in \cite{LL} is equivalent to Definition~\ref{def:RRFscalar}.

We now discuss some open questions in the pointwise setting. One feature of Gromov's definition is that it may be localized around a~point $x\in M$, by requiring only that there exist a~neighborhood of~$x$ such that no cube within the neighborhood that contains $x$ has strictly mean convex faces and acute dihedral angles; see \cite[p.~1144]{Gro}. In light of this, for $x\in M$ define~$C^0_{\G}(x, \kappa)$ to be the space of germs of~$C^0$ metrics on $M$ at $x$ that have scalar curvature bounded below by~$\kappa$ at~$x$ in the sense of~\cite{Gro}.

Define $C^0_{\beta}(x,\kappa)$ to be the space of germs of $C^0$ metrics on $M$ at $x$ that have scalar curvature bounded below by $\kappa$ at $x$ in the sense of Definition~\ref{def:RRFscalar}. By Theorem~\ref{thm:preservation}, this is a reasonable of localization of $C^0_{\Rf}(M,\kappa)$, as was intended.

\begin{Question}
Let $M$ be a closed manifold and $x\in M$. Do we have $C^0_{\G}(x, \kappa) = C^0_{\beta}(x,\kappa)$?
\end{Question}
\begin{Remark}
In contrast with the case of global lower bounds, the results in this paper do \emph{not} imply that $C^0_{\beta}(x,\kappa)\subset C^0_{\G}(x,\kappa)$.
\end{Remark}

Theorem \ref{thm:secondorderagreement} says that the perturbation of a metric by greater than second order does not affect $\beta$-weak lower bounds on the scalar curvature. Thus, it is natural to ask whether one can characterize, up to higher order perturbation, those metrics, on $\R^n$, say, that have nonnegative scalar curvature in the sense of Definition \ref{def:RRFscalar}:

\begin{Question}
Suppose $g = g_{ij}{\rm d}x^i\otimes {\rm d}x^j$ is a metric on a neighborhood of the origin in $\R^n$, and that we may write
\begin{equation*}
g_{ij}(x) = \delta_{ij} + r^2G_{ij}\left(\frac{x}{r}\right) + O\big(|x|^{2+\eta}\big),
\end{equation*}
where the $G_{ij}$ are functions on $\mathbb{S}^{n-1}\subset \R^n$ satisfying $x^ix^jG_{ij}(x) = 0$. Is there an explicit characterization of metrics of this form that have nonnegative scalar curvature at the origin, in the sense of Definition~{\rm \ref{def:RRFscalar}}?
\end{Question}

\subsection*{Acknowledgements}
This survey is dedicated to Misha Gromov on the occasion of his 75th birthday. I would like to thank him for his interest in my work. I would also like to thank my advisor, Richard Bamler, for his guidance during the writing process. Finally, I would like to thank the referees for numerous helpful comments on a previous draft of this paper.
This material is based upon work supported by the National Science Foundation Graduate Research Fellowship Program under Grant No.~DGE 1752814. %Any opinions, findings, and conclusions or recommendations expressed in this material are those of the author(s) and do not necessarily reflect the views of the National Science Foundation.

\pdfbookmark[1]{References}{ref}
\LastPageEnding

\end{document}